\documentclass[11pt]{article}
\usepackage{amssymb}
\usepackage{amsmath}
\usepackage{amsthm}
\usepackage{epsfig}
\usepackage{mathrsfs}

\numberwithin{equation}{section}

\newtheorem{thm}{Theorem}[section]
\newtheorem{lem}[thm]{Lemma}

\newtheorem{cor}[thm]{Corollary}

\newtheorem{rem}[thm]{Remark}
\newtheorem{exa}{Example}[section]

\newcommand{\be}{\begin{equation}}
\newcommand{\ee}{\end{equation}}
\newcommand{\bea}{\begin{eqnarray*}}
\newcommand{\eea}{\end{eqnarray*}}
\newcommand{\ben}{\begin{eqnarray}}
\newcommand{\een}{\end{eqnarray}}
\newcommand{\R}{\mathbb{R}^2}
\newcommand{\Z}{\mathbb{Z}^2}
\newcommand{\tp}{\tilde{\psi}}

\newcommand{\ta}{\tilde{A}}

\newcommand{\mR}{\mathbb{R}}

\newcommand{\ml}{\mathcal{L}}

\newcommand{\mc}{\mathcal{C}}

\newenvironment{keyword}
{\noindent{\textbf{Keywords.}}}

\newenvironment{class}
{\noindent{\textbf{AMS subject classification:}}}

\newenvironment{ac}{\noindent{\bf Acknowledgements}}

\begin{document}
\date{}
\title{Shearlet frames with short support}
\author{
 Song Li\thanks{Department of Mathematics,
 Zhejiang University Hangzhou, 310027, China}
 and
 Yi Shen\thanks{Corresponding author: sy1133@163.com, Department of Mathematics,
 Zhejiang University Hangzhou, 310027, China}}
\maketitle

\begin{abstract}
Compactly supported shearlets have been studied in both theory and
applications. In this paper,  we construct symmetric compactly
supported shearlet systems based on pseudo splines of type II.
Specially, using B-splines, we construct shearlet frame having
explicit analytical forms which is important for applications. The
shearlet systems based on B-splines also provide optimally sparse
approximation within cartoon-liked image.

\bigskip
\bigskip
\begin{class}
Prime 42C40; Secondary: 42C15, 65T60, 65T99, 94A08
\end{class}
\end{abstract}

\section{Introduction}
Cartoon-liked image are 2-dimensional functions that are $C^2$
except for discontinuities along $C^2$ curves \cite{CD}. To find
optimally sparse representations of cartoon-like image, several
variations of the wavelet scheme have been  proposed, such as
curvelets \cite{CD} and contourlets \cite{DV}. Shearlets frame
developed in \cite{LLKW}  is  the first multiscale directional
system which also provides almost optimally sparse approximation
with cartoon-like images. However, these studies are only concerned
band limited generator. Very recently, Kutyniok and Lim presented a
complete proof of (almost) optimally sparse approximations of
cartoon liked images by using shearlet systems which are generated
by compactly supported shearlets under some weak moment conditions
\cite{KL}. They also constructed a class of compactly supported
shearlet frames based on pseudo splines of type I \cite{KGL2}.
Hence, excellent spatial localization is achieved. But the shearlet
frame still have two disadvantages:
\begin{itemize}
\item The shearlet  is not symmetric or anti-symmetric;
\item The shearlets do not have explicit analytical forms in spatial
domain.
\end{itemize}
These drawbacks motivate us to consider constructing shearlet frame
using B-splines. B-splines had a significant impact on the
development of the theory of the wavelet analysis. They yield the
only wavelets that have explicit analytical forms. All other wavelet
bases are defined indirectly through an infinite product in Fourier
domain \cite{C,DTen}. As Daubechies pointed out in \cite{DTen},
except the Haar wavelet function,  there is no compactly supported
real-valued symmetric orthonormal  wavelet basis in $L_2(\mathbb
R)$. However, it is much easier and more flexible to construct and
design compactly supported wavelet frames or Riesz bases than
orthonormal wavelet bases. For example, from any B-spline function
of order $m$, one can construct a symmetric tight wavelet frame with
$m$ generators \cite{RS}. Tight wavelet frame from B-splines with
high vanishing moments were considered in \cite{CH,DHRS}. The
compactly supported Riesz wavelets generated from B-splines were
first constructed in \cite{CW}. The shortest supported Riesz wavelet
with $m$ vanishing moments from B-spline of order $m$ were
constructed in \cite{HSshort}. The compactly supported wavelet bases
from B-splines for Sobolev spaces were investigated in
\cite{HSCA,JWZ}.

The rest of the paper is organized as follows. In Section
\ref{S:Bsplines}, we construct shearlets  based on B-splines, then
we present some results on the optimally sparse approximations of
cartoon-like images. In Section 3, we investigate the lower bounds
and the upper bounds of the pseudo splines in Fourier domain. These
results not only have their own interests but also have  closed
relations to the shearlets frame bounds. In sections 4, we presents
some examples to illustrate our results.

\section{Shearelts based on B-splines}\label{S:Bsplines}

In this section, we  first introduce the definitions of shearlet
frame and B-splines function. Then we construct compactly supported
shearlet frame with the generator from B-splines. Finally, we show
that these compactly supported shearlet systems provide (almost)
optimally sparse approximations of cartoon-liked images.

Shearlets are scaled according to a parabolic scaling law encoded in
the matrix $A_{2j}$ or $\tilde{A}_{2j}$, $j\in\mathbb Z$, and
exhibit directionality by parameterizing slope encoded in the
matrices $S_k$, $k\in\mathbb Z$, defined by
$$
A_{2^j}=\left( \begin{array}{cc}
 2^j& 0 \\
 0 & 2^{j/2}
\end{array}
\right),
 \quad
\tilde{A}_{2^j}=\left( \begin{array}{cc}
 2^{j/2}& 0 \\
 0 & 2^{j}
\end{array}
\right),
 \quad
 {S}_k=
 \left(
 \begin{array}{cc}
 1 & k\\
 0 & 1
 \end{array}
 \right).
$$
Now  we define discrete shearlet systems in $2$D.  Let
$c=(c_1,c_2)\in (\mathbb R_+)^2$. For $\phi, \psi, \tilde{\psi}\in
L^2(\R)$, the cone-adapted $2$D discrete shearlet system
  $SH(\phi,\psi,\tilde{\psi};c)$ is defined by
  $$
  SH(\phi,\psi,\tilde{\psi};c) = \Phi(\phi;c_1)\cup \Psi(\psi;c) \cup
  \tilde{\Psi}(\tp;c),
  $$
 where
  \bea
  \Phi(\phi;c_1) &=& \{\phi(\cdot-cm): m\in \Z\},\\
  \Psi(\psi;c) &=& \{ 2^{\frac{3}{4}j} \psi (S_k A_{2^j}\cdot-cm): \
  j\geq 0, \  | k| \leq \lceil 2^{j/2} \rceil, \  m\in\Z
  \},\\
  \tilde{\Psi}(\tilde{\psi};c) &=& \{ 2^{\frac{3}{4}j} \tilde{\psi} (S^T_k \ta_{2^j}\cdot-cm):
  \  j\geq 0, \  |k| \leq  \lceil 2^{j/2} \rceil,
  \   m\in\Z
  \}.
  \eea
We partite the frequency plane into
$\mathcal{C}_1(\alpha)-\mathcal{C}_4(\alpha)$ where
$$
\mathcal{C}_{l}(\alpha) = \left\{
\begin{array}{ll}
\{(\xi_1,\xi_2)\in\R:\xi_1 \geq \alpha,\ |\xi_2/\xi_1|\leq 1 \} :  & l=1\\
\{(\xi_1,\xi_2)\in\R:\xi_2 \geq \alpha,\ |\xi_1/\xi_2|\leq 1 \} :  & l=2\\
\{(\xi_1,\xi_2)\in\R:\xi_1 \geq -\alpha,\ |\xi_2/\xi_1|\leq 1 \} :  & l=3\\
\{(\xi_1,\xi_2)\in\R:\xi_2 \geq -\alpha,\ |\xi_1/\xi_2|\leq 1 \} :  & l=4\\
\end{array}
\right.
$$
and a centered rectangle
$$
\mathcal{R}(\alpha) = \{(\xi_1,\xi_2)\mathbb{R}^2:
\|(\xi_1,\xi_2)\|_{\infty}<\alpha\}.
$$
The region $\mc_{1}\cup \mc_{3}$ is covered by the frequency support
of shearelets in $\Psi(\psi;c)$. The region $\mc_{2}\cup \mc_{4}$ is
covered by the frequency support of shearelets in
$\tilde{\Psi}(\tilde{\psi};c)$. The region $\mathcal{R}$ is covered
by the frequency support of $\phi$. Recall  $\{\sigma_i\}_{i\in I}$
form a frame in $L^2(\R)$ if there exist constants $A, B
> 0$ such that
 $$
 A\|f\|_2^2\leq \sum_{i\in I} |\langle f, \sigma_i \rangle|^2
 \leq
 B\|f\|_2^2 \quad \forall\ f\in L^2(\R).
 $$
The numbers $A$, $B$ are called frame bounds. If
$SH(\phi,\psi,\tilde{\psi};c)$ is a frame for $L^2(\R)$, we refer to
$\psi$ and $\tp$ as \emph{shearlets}.

We say that $\phi$ is a \emph{refinable} function with \emph{mask}
$\hat{a}(\xi)$ if $ \hat{\phi}(2\xi) = \hat{a}(\xi)\hat{\phi}(\xi)$,
$\xi\in\mathbb R^d$. The \emph{Fourier transform} of a function
$f\in L_1(\mathbb R^d)$ is defined as $ \hat{f}(\xi)=
 \int_{\mathbb{R}^d}f(x)e^{-ix\cdot\xi}dx$
and can be naturally extended to tempered distributions. As an
important family of refinable functions, B-spline functions are
useful in applications.  B-spline with order $m$ and its mask is
defined by
$$
\widehat{B_m}(\xi)=
e^{\frac{-ij\xi}{2}}\biggl(\frac{\sin(\xi/2)}{\xi/2}\biggr)^m\quad
\mathrm{and}\quad \widehat{a}(\xi)=
e^{\frac{-ij\xi}{2}}\cos^m(\xi/2), \quad \xi\in \mathbb{R},
$$
where $j=0$ when $m$ is even, and $j=1$ when $m$ is odd. The
B-spline function $B_m\in C^{m-2}(\mathbb R)$ is a function of
piecewise polynomials of degree less than $m$, vanishes outside the
interval $[0,m]$ and is symmetric about the point $x=m/2$. Now we
state our first contribution  on the shearlet frame in this paper.

\begin{thm}\label{Bframe}
Let $B_N$ be the B-spline function of order $N$ with the mask $
\widehat{a_N}(\xi)= (2^{-N})(1+e^{-i\xi})^N$. Let $N_1,N_2\in
\mathbb {N}$ be such that $N_1>N_2>3$. Define a shearlet $\psi\in
L^2(\mathbb R^2)$ by
$$
\hat{\psi}(2\xi) = 2^{-N_1}e^{-i\xi_1} (1-e^{-i\xi_1})^{N_1}
\widehat{B_{N_2}}(\xi_1)\widehat{B_{N_2}}(\xi_2),\quad
\xi=(\xi_1,\xi_2).
$$
For given $0<\alpha<\pi/2$, there exits a sampling constant
$\hat{c}>0$ such that the shearlet system $\Psi(\psi;c)$ forms a
frame for
 $
 \{ f\in L_2(\R) :
 \mathrm{supp} \hat{f}\in \mathcal{C}_1(\alpha)\cup\mathcal{C}_3(\alpha)  \}
 $
for  $c_2\leq c_1\leq \hat{c}$.
\end{thm}
\begin{proof}
The proof is a straightforward consequence of Theorem
\ref{thmframe}.
\end{proof}

The cartoon-liked model was first introduced in \cite{CD}. The basic
idea is to choose a closed boundary curve and then fill the interior
and exterior parts with $C^2$ functions. For $\nu>0$, the set
$STAR^2(\nu)$  is defined to be the set of all $B\subset[0,1]^2$
such that $B$ is a translate of a set
$$
\{ x\in \mathbb{R}^2:  |x|\leq \rho(\theta),\ x = (|x|,\theta) \
\text{in polar coordinates} \}
$$
which satisfies $|\rho''(\theta)<\nu|,\ \rho\leq \rho_0<1$. Then
$\mathscr{E}^2(\nu)$ denotes the set of functions $f\in L(\R)$ of
the form
$$
f=f_0+f_1\chi_{B},
$$
where $f_0,$ $f_1\in C_0^2([0,1]^2)$
and $B\in STAR^2(\nu)$. The bandlimited curvelets, contourlets, and
shearlets exhibit (almost) optimally sparse approximation with this
model. The first complete proof of (almost) optimally sparse
approximations of cartoon-liked images by using compactly supported
shearlet frame was given in \cite{KL}. Let us now be more precise,
and introduce these results. Let $c>0$, and let $\phi$, $\psi$,
$\tilde{\psi}\in L_2(\R)$ be compactly supported. Suppose that for
all $\xi=(\xi_1,\xi_2)\in \R$, the shearlet $\psi$ satisfies
    \be\label{condition1}
    |\hat{\psi}(\xi)|\leq C_1\cdot \min\{1,|\xi_1|^{\alpha}\}\cdot \min\{1,|\xi_1|^{-\gamma}\}
    \cdot \min\{1,|\xi_2|^{-\gamma}\},
    \ee
    \be\label{condition2}
    |\frac{\partial }{\partial \xi_2}\hat{\psi}(\xi)|\leq
    |h(\xi_1)|\left(1+\frac{|\xi_2|}{|\xi_1|}\right)^{-\gamma},
    \ee
where $\alpha>5$, $\gamma\geq 4$, $h\in L_1(\mathbb R)$, and $C_1$
is a constant, and suppose that the shearlet $\tp$ satisfies
(\ref{condition1}) and (\ref{condition2}) with the roles of $\xi_1$
and $\xi_2$ reversed. Further, suppose that $SH(c;\phi,\psi,\tp)$
forms a frame for $L_2(\R)$. Denote $(\sigma_i)_{i\in I}=
SH(c;\phi,\psi,\tp)$. Let $(\tilde{\sigma})_{i\in I}$ be a dual
frame of $(\sigma_i)_{i\in I}$. We can take the nonlinear $N$-terms
approximation
$$
f_N = \sum_{i\in I_N} \langle f, \sigma_i\rangle \tilde{\sigma}_i,
$$
where $(\langle f, \sigma_i\rangle)_{i\in I_N}$ are the $N$ largest
coefficients $\langle f, \sigma_i\rangle$ in magnitude. Then, for
any $v>0$, the shearlet frame $SH(c;\phi,\psi,\tp)$ provides
(almost) optimal sparse approximation of function
$f\in\mathscr{E}^2(v)$ in the sense that there exists some $C>0$
such that
  $$
  \|f-f_N\|^2_2 \leq C \cdot(\log N)^3 \cdot N^{-2}.
  $$
The condition (\ref{condition1}) and (\ref{condition2}) can be
viewed as a generalization of a second order directional vanishing
moment condition, which is crucial for having fast decay of the
shearlet coefficients. Following the line of \cite{KGL2}, we obtain
the following results which is our second contribution in this
paper.
\begin{thm}
Let $B_N$ be the B-spline function of order $N$ with the mask $
\widehat{a_N}(\xi)= (2^{-N})(1+e^{-i\xi})^N$. Let $N_1,N_2\in
\mathbb {N}$ be such that $N_1>5$, $N_2\geq 4$ and $N_1>N_2$. Define
a shearlet $\psi_1\in L^2(\mathbb R^2)$ by
$$
\widehat{\psi_1}(2\xi) = 2^{-N_1}e^{-i\xi_1} (1-e^{-i\xi_1})^{N_1}
\widehat{B_{N_2}}(\xi_1)\widehat{B_{N_2}}(\xi_2)\quad
\xi=(\xi_1,\xi_2).
$$
Let $\phi(x) = B_{N_2}(x_1)B_{N_2}(x_2)$ and
$\psi_2(x_1,x_2)=\psi_1(x_1,x_2)$. Then there exit  sampling
constant $c>0$ such that the shearlet system $SH(\phi,
\psi_1,\psi_2;c)$ provides (almost) optimally sparse approximations
of function $f\in\mathscr{E}^2(v)$ in the sense that there exists
some constant $C>0$ such that
  $$
  \|f-f_N\|^2_2 \leq C \cdot(\log N)^3 \cdot N^{-2},
  $$
where $f_N$ is the nonlinear $N$-term approximation obtained by
choosing the $N$ largest shearlet coefficients of $f$.
\end{thm}
\begin{proof}
By Theorem \ref{Bframe}, there exits a sampling constant $\hat{c}>0$
such that the shearlet system $\Psi(\psi;c)$ forms a frame for
 $
 \{ f\in L_2(\R) :
 \mathrm{supp} \hat{f}\in \mathcal{C}_1(\alpha)\cup\mathcal{C}_3(\alpha)  \}
 $
 with $c=(c_1,c_2)\in (\mathbb R_+)^2$
and $c_2\leq c_1\leq \hat{c}$. With the same argument as in
\cite{KGL2}, we can prove that $SH(c;\phi,\psi,\tp)$ forms a frame
for $L_2(\mR)$ with  $c=(c_1,c_2)\in (\mathbb R^+)^2$ and $c_2\leq
c_1\leq \hat{c}$.

Rewrite $\widehat{\psi_1}(\xi)$ in the following
$$
\widehat{\psi_1}(\xi) = 2^{-N_1}e^{-i\xi_1}
(1-e^{-i\xi_1})^{N_1-N_2}\widehat{B_{N_2}}(\xi_1)\cdot
(1-e^{-i\xi_1})^{N_2} \widehat{B_{N_2}}(\xi_2).
$$
We obtain
 \bea
\left|\frac{\partial\widehat{\psi_1}(\xi)}{\partial \xi_2}\right| &
= & \left|2^{-N_1}e^{-i\xi_1}
(1-e^{-i\xi_1})^{N_1-N_2}\widehat{B_{N_2}}(\xi_1)\right|\cdot
\left|(1-e^{-i\xi_1})^{N_2} \widehat{B'_{N_2}}(\xi_2)\right| \\
 &\leq &
 \left|2^{-N_1}e^{-i\xi_1}
(1-e^{-i\xi_1})^{N_1-N_2}\widehat{B_{N_2}}(\xi_1)\right|\cdot
|\xi_1|^{N_2} \cdot \left|\widehat{B'_{N_2}}(\xi_2)\right|.
 \eea

Set $h(\xi_1) = 2^{-N_1}e^{-i\xi_1}
(1-e^{-i\xi_1})^{N_1-N_2}\widehat{B_{N_2}}(\xi_1)$. We have
$$
|h(\xi_1)| \leq C'' \widehat{B_{N_2}}(\xi_1) \leq C''
|\xi_1|^{-N_2}.
$$
Therefore, $h(\xi_1)\in L_1(\mathbb R)$.

Let us consider $|\widehat{B_{N_2}}'(\xi_2)|$. When $N_2$ is even, $
\widehat{B_{N_2}}(\xi_2) = \left(  \frac{\sin(\xi_2/2)}{\xi_2/2}
\right)^{N_2}$. Then
 \bea
\widehat{B_{N_2}}'(\xi_2) &=& N_2 \left(
\frac{\sin(\xi_2/2)}{\xi_2/2} \right)^{N_2-1}\left[
\frac{\cos(\xi_2/2)}{4\xi_2} - \frac{\sin(\xi_2/2)}{2\xi_2^2}
\right]\\
&\leq& C |\xi_2|^{-N_2}.
 \eea
When $N_2$ is odd, $ \widehat{B_{N_2}}(\xi_2) = e^{-i\xi_2}\left(
\frac{\sin(\xi_2/2)}{\xi_2/2} \right)^{N_2}$. Then
 \bea
 |\widehat{B_{N_2}}'(\xi_2)|& =&
 \left|-ie^{-i\xi_2}\left(
 \frac{\sin(\xi_2/2)}{\xi_2/2} \right)^{N_2}
 +e^{-i\xi_2}N_2 \left(
 \frac{\sin(\xi_2/2)}{\xi_2/2} \right)^{N_2-1}\left[
 \frac{\cos(\xi_2/2)}{4\xi_2} - \frac{\sin(\xi_2/2)}{2\xi_2^2}
 \right] \right| \\
 &\leq&
 \left|  \left(
 \frac{\sin(\xi_2/2)}{\xi_2/2} \right)^{N_2}  \right| +
 \left|  N_2 \left(
 \frac{\sin(\xi_2/2)}{\xi_2/2} \right)^{N_2-1}\left[
 \frac{\cos(\xi_2/2)}{4\xi_2} - \frac{\sin(\xi_2/2)}{2\xi_2^2}
 \right]   \right| \\
 &\leq& C' |\xi_2|^{-N_2}.
 \eea
We can conclude that   there exists a constant $C$ such that
$$
|\widehat{\psi_1}(\xi)|\leq C \min\{1, |\xi_1|^{N_1} \}\cdot
\min\{1, |\xi_1|^{-N_2} \}\cdot \min\{1, |\xi_2|^{-N_2}\},
$$
and
$$
\left|\frac{\partial\widehat{\psi_1}(\xi)}{\partial \xi_2}\right|
\leq |h(\xi_1)|\left( 1 + \frac{|\xi_2|}{|\xi_1|} \right)^{-N_2},
$$
where  $N_1>5$, $N_2 \geq 4$, $h\in L_1(\mathbb R)$. Hence, the
shearlet $\psi_1$ satisfies the conditions (\ref{condition1}) and
(\ref{condition2}), and the shearlet $\psi_2$ likewise. The theorem
is proved
\end{proof}

\section{Shearlets based on pseudo splines}
In this section, we first briefly recall a family of refinable
function: pseudo splines.  After establishing some useful lemmas, we
 investigate  the lower bound and  the decay of the Fourier
transform of pseudo splines of type II. Finally, we construct
shearlet frame based on pseudo splines.

Pseudo-splines are compactly supported refinable functions in
$L_2(\mR)$. For positive integers $N,l\in\mathbb N$ with $l<N$,
denote $ P_{N,l}(x):=\sum_{j=0}^l\binom{N-1+j}{j}x^j$. The  mask of
a pseudo spline of type II with order $(N,l)$ is defined by
 \be\label{mask}
 \widehat{_2a_{N,l}}(\xi)
 :=\cos^{2N}(\xi/2)P_{N,l}(\sin^2(\xi/2)).
 \ee
The mask of a pseudo spline of type I is defined by
 $\widehat{_2a_{N,l}}(\xi):=|\widehat{_1a_{N,l}}(\xi)|^2$.
Hence, $\widehat{_1a_{N,l}}(\xi)$ with real coefficients is the
square root of $|\widehat{_1a_{N,l}}(\xi)|^2$  using the Riesz
lemma. In general, pseudo splines of type I are neither symmetric
nor antisymmetric. To achieve symmetry,  compactly supported complex
valued pseudo splines were introduced in \cite{SLM}. The
corresponding pseudo splines can be defined in terms of their
Fourier transform, i.e.
$$
\widehat{_k\phi}(\xi) = \prod_{j=1}^{\infty}
\widehat{_ka_{N,l}}(2^{-j}\xi),\quad k=1,2.
$$
An important fact is that  $\widehat{_2a_{N,l}}(\xi)$ is defined by
the summation of the first $l+1$ terms of the binomial expansion of
$ (\cos^2(\xi/2)+\sin^2(\xi/2))^{N+l}=1$ \cite{DHRS,DS1} i.e.
 \be\label{mask_a_equal}
\sum_{j=0}^l\binom{N-1+j}{j}\sin^{2j}(\xi/2)=
\sum_{j=0}^l\binom{N+l}{j}\sin^{2j}(\xi/2)\cos^{2(l-j)}(\xi/2).
 \ee
The pseudo splines with order $(N, 0)$  are B-splines. For the case
$l=N-1$, the pseudo splines of type I are orthogonal refinable
functions  given in \cite{DTen}, and the pseudo splines of type II
are interpolatory refinable function given in \cite{Dubuc}. The
other pseudo splines fill in the gap between the B-spline and
orthogonal refinable functions for type I and B-spline and
interpolatory refinable function for type II.
For positive integers $N,l\in\mathbb N$  with $l<N$, let
$\hat{a}(\xi)$ be the mask of the pseudo splines of type II with
order $(N,l)$. Then $\widehat{a_{N,l}}(\xi)$ can be factorized as
 \be\label{Lfactor}
 |\widehat{a_{N,l}}(\xi)| = \cos^{2N}(\xi/2)|\mathcal{L}(\xi)|,\quad
 \xi\in[-\pi,\pi].
 \ee
This shows that pseudo splines is the convolution of a B-spline of
some order with a distribution. Since $\ml(\xi)$ is bounded,
$\ml(\xi)$ is actually the mask of a refinable distribution. The
regularity of $\phi$ comes from the $\cos^{2N}(\xi/2)$ factor. The
distribution part provides some desirable properties for $\phi$,
such as  orthogonality of its shifts. In \cite{DS1}, Dong and Shen
gave a regularity analysis of pseudo splines of both types. The key
to regularity analysis is
 \ben\label{l1}
& &|\mathcal{L}(\xi)| \leq \left|\mathcal{L}(\frac{2\pi}{3})\right|,
\qquad
|\xi|\leq \frac{2\pi}{3},\\
& & |\mathcal{L}(\xi)\mathcal{L}(2\xi)|    \leq
|\mathcal{L}(\frac{2\pi}{3})|^{2}, \qquad  \frac{2\pi}{3} \leq
|\xi|\leq \pi\nonumber.
 \een
To investigate the bounds and the decay of the pseudo splines in
Fourier domain, we establish the following lemmas which have their
own interesting. For simplify, we   denote the mask of pseudo
splines of type II by $\hat{a}(\xi)$.
 \subsection{Lemmas}
\begin{lem}
For positive integers $N,l\in\mathbb N$  with $l<N$, let
$\hat{a}(\xi)$ be the mask of the pseudo splines of type II with
order $(N,l)$ and let $\ml(\xi)$ be defined as in (\ref{Lfactor}).
Then
 $$
 |\hat{a}(\xi)| \geq 1- C_1 |\xi|^{2l+2},
 $$
where $C_1=\frac{\sum_{j=l+1}^{N+l}\binom{N+l}{j}}{2^{2l+2}}$.
\end{lem}
\begin{proof}
By (\ref{mask_a_equal}), we obtain
 \bea
 1-|\hat{a}(\xi)| & = &\cos^{2N}(\xi/2) \sum_{j=l+1}^{N+l}
 \binom{N+l}{j}\sin^{2j}(\xi/2)\cos^{2(l-j)}(\xi/2)\\
 &\leq& \sin^{2l+2}(\xi/2) \sum_{j=l+1}^{N+l}
 \binom{N+l}{j} \\
 &\leq&
 \frac{\sum_{j=l+1}^{N+l}\binom{N+l}{j}}{2^{2l+2}}|\xi|^{2l+2}.
 \eea
\end{proof}
\begin{lem}
For positive integers $N,l\in\mathbb N$  with $l<N$, let
$\hat{a}(\xi)$ be the mask of the pseudo splines of type II with
order $(N,l)$. Then
 $$
 |\mathcal{L}(\xi)| \leq 1 + C_2 |\xi|^{2},
 $$
where $C_2=\frac{\sum_{j=1}^{l}\binom{N-1+j}{j}}{4}$.
\end{lem}
\begin{proof}
We obtain
  \bea
  |\ml(\xi)| -1
  & = &  \sum_{j=0}^{l}\binom{N-1+j}{j} \sin^{2l}(\xi/2)-1 \\
  & = & \sin^{2}(\xi/2)\sum_{j=1}^{l}\binom{N-1+j}{j}
  \sin^{2l-2}(\xi/2)\\
  &\leq& \frac{\sum_{j=1}^{l}\binom{N-1+j}{j}}{4}|\xi|^2.
  \eea
\end{proof}

\begin{lem}\label{lem_b_boud}
For positive integers $N,l\in\mathbb N$  with $l<N$, let
$\hat{a}(\xi)$ be the mask of the pseudo splines of type II with
order $(N,l)$. Let $\hat{b}(\xi) =
e^{-i\xi}\overline{\hat{a}(\xi+\pi)}$, then
 \be\label{blowup}
 |\hat{b}(\alpha)| \chi_{[-\beta, -\alpha]\cup[\alpha,\beta]}(\xi)
 \leq|\hat{b}(\xi)| \leq
 \min\{1, \frac{\sum_{j=0}^l\binom{N+l}{j}}{2^{2N}}|\xi|^{2N}\},\quad 0\leq
 \alpha< \beta\leq \pi.
 \ee
\end{lem}
\begin{proof}
$|\hat{b}(\xi)|$ is a $\pi$ shift of the $2\pi$-periodic triangle
polynomial $\hat{a}(\xi)$. Note that $|\hat{a}(\xi)|$ is increasing
on $[-\pi,0]$ and decreasing on $[0,\pi]$. It is easy to see that
$|\hat{b}(\xi)|$ is decreasing on $[-\pi,0]$ and increasing on
$[0,\pi]$. Hence the left side of (\ref{blowup}) holds. By
(\ref{Lfactor}),  We have
  \bea
  |\hat{b}(\xi)| &=& \sin^{2N}(\xi/2) \sum_{j=0}^{l}
  \binom{N+l}{j}\cos^{2j}(\xi/2)\sin^{2(l-j)}(\xi/2)\\
  &\leq& \frac{\sum_{j=0}^l\binom{N+l}{j}}{2^{2N}}|\xi|^{2N}.
  \eea
\end{proof}

\subsection{Regularity}
Now we give the lower bound for the pseudo splines of type II in the
Fourier domain.
\begin{thm}\label{thmphilow}
  Let $\phi$ be the pseudo-splines of Type II with order $(N,l)$.
  Let $K \subset [-\pi,\pi] $, then
  $$
  |\hat{\phi}(\xi)| \geq C_4 \cdot \chi_{K}(\xi).
  $$
where
$C_4=\prod_{k=1}^{k_0}|\hat{a}(2^{-k}\xi_0)|\exp(-C_12^{-k_0+1}|\xi_0|^{2l+2})$.
\end{thm}
\begin{proof}
Let $\xi_0 = \arg\max_{\xi\in K}|\xi|$. Since $\hat{a}(\xi)$ is
decreasing on $[0, \pi]$, we have
 $$
 |\hat{a}(2^{-k}\xi)| \geq   |\hat{a}(2^{-k}\xi_0)|
 $$
for $k\geq 1$ and $\xi\in K$. We choose sufficiently large $k_0$ so
that $2^{-k}C_1|\xi|^{2l+2}<\frac{1}{2}$ if $\xi\in K$ and $k>k_0$.
Using $1-x \geq e^{-2x}$ for $0\leq x\leq 1/2$, we obtain
 \bea
 |\hat{a}(\xi)| \geq 1 - C_1 2^{-k} |\xi|^{2l+2} \geq \exp(-2^{k+1}C_1|\xi|^{2l+2}),\quad k>k_0.
 \eea
Therefore, for $\xi\in K$
 \bea
 |\hat{\phi}(\xi)|& =& \prod_{k=1}^{k_0}|\hat{a}(2^{-k}\xi)|
                       \prod_{k=k_0+1}^{\infty}|\hat{a}(2^{-k}\xi)|\\
 &\geq& \prod_{k=1}^{k_0}|\hat{a}(2^{-k}\xi_0)|\prod_{k=k_0+1}^{\infty}
 \exp(-2C_12^{-k}|\xi|^{2l+2})\\
 &\geq&\prod_{k=1}^{k_0}|\hat{a}(2^{-k}\xi_0)|\exp(-C_12^{-k_0+1}|\xi_0|^{2l+2})=C_4.
 \eea
\end{proof}

In the following, we mainly investigate the decay of  the pseudo
splines of type II in the Fourier domain.
\begin{lem}
For positive integers $N,l\in\mathbb N$  with $l<N$, let
$\hat{a}(\xi)$ be the mask of the pseudo splines of type II with
order $(N,l)$. Let $\ml(\xi)$ be the distribution part. Assume that
$\widehat{\phi_{\ml}}(\xi)=\prod_{j=1}^{\infty}\ml(2^{-j}\xi)$, then
 \be\label{decayL}
 \widehat{\phi_{\ml}}(\xi)\leq \exp(C_2/3) q_1q_2^{J-1} |\xi|^{\log_2(q_1^{1/(J-1)}q_2)},
 \ee
where $q_1 = \sup_{|\xi|\leq
\pi}|\ml(\xi)|=\sum_{j=0}^l\binom{N-1+j}{j}$ and $q_2 =
|\ml(2\pi/3)|$.
\end{lem}
\begin{proof}
  Since $|\mathcal{L}(\xi)| \leq 1 + C_2 |\xi|^{2}$, we have for
$|\xi|\leq 1$,
 \bea
 \prod_{j=1}^{\infty}\left| \ml (2^{-j}\xi) \right|
 &=& \exp\left(\sum_{j=1}^{\infty} \log (| \ml (2^{-j} \xi| )\right) \\
 &\leq& \exp\left(C_2 \sum_{j=1}^{\infty} 2^{-2j}|\xi| \right)\\
 &\leq& \exp(C_2/3).
 \eea
Now for $|\xi| >1$, we have
 $$
 \prod_{k=1}^{\infty} \ml(2^{-j}\xi) =  \prod_{k=0}^{\infty}f(2^{-Jk}\xi)
 $$
where $ f(\xi) =\prod_{j=1}^J\ml(2^{-j}\xi)$. For the given $|\xi|
>1$, there exists a positive integer $k_0$ such that $2^{k_0J} \leq
|\xi| \leq 2^{(k_0+1)J}$. Denote $\eta = 2^{-(k_0+1)J}\xi$, we
obtain $|\eta|<1$. Define $q_1 = \sup_{|\xi|\leq \pi}|\ml(\xi)|$ and
$q_2 = |\ml(2\pi/3)|$. By (\ref{l1}), we obtain
$$
\left|\prod_{j=1}^J\ml(2^{-j}\xi)\right| \leq q_1 q_2^{J-1}.
$$
Therefore,
 \bea
 \prod_{k=k_0+1}^{\infty} |f(2^{-kJ}\xi)| & = &
 \prod_{j=0}^{\infty} |f(2^{-(j+k_0+1)J}\xi)|
 = \prod_{j=0}^{\infty} |f(2^{-jJ}2^{-(k_0+1)J}\xi)|
 = \prod_{j=0}^{\infty} \left| f(2^{-jk}\eta) \right|  \\
 &=& \prod_{j=0}^{\infty} \left| \ml(2^{-j}\eta) \right| \leq \exp(C_2/3).
 \eea
Moreover,
 \bea
 \prod_{k=0}^{k_0}|f(2^{-kJ})|
 &\leq &
 \left(q_1q_2^{J-1}\right)^{k_0+1}
 \leq
 q_1q_2^{J-1}\left(q_1^{1/(J-1)}q_2\right)^{k_0(J-1)}\\
 &\leq&
 q_1q_2^{J-1}\left(q_1^{1/(J-1)}q_2\right)^{\log_2|\xi|}\\
 &\leq&
 q_1q_2^{J-1} |\xi|^{\log_2(q_1^{1/(J-1)}q_2)}.
 \eea
We conclude that
$$
\prod_{j=1}^{\infty} \ml(2^{-j}\xi) = \exp(C_2/3)  q_1q_2^{J-1}
|\xi|^{\log_2(q_1^{1/(J-1)}q_2)}.
$$
\end{proof}

\begin{thm}\label{phiup}
For positive integers $N,l\in\mathbb N$  with $l<N$, let $\phi$ be
the pseudo splines of type II with order $(N,l)$. Then for any given
integer $J\geq 1$.
  \be\label{phic}
  |\hat{\phi}(\xi)| \leq \min\{1, C_3
  |\xi|^{-2N+\log_2(q_1^{1/(J-1)}q_2)}\}
  \ee
where $C_3 = 4^N\exp(C_2/3) q_1q_2^{J-1}$.
\end{thm}

\begin{proof}
Since
 $$
 |\hat{a}(\xi)|^2+|\hat{a}(\xi+\pi)|^2 \leq
 |\hat{a}(\xi)|+|\hat{a}(\xi+\pi)| \leq
 1,
 $$
we obtain $|\hat{\phi}(\xi)|\leq 1$.

It is well known that
$$
\prod_{j=1}^{\infty}\cos(2^{-j}\xi) = \frac{\sin(\xi/2)}{\xi/2}.
$$
Hence
$$
|\hat{\phi}(\xi)| = \left(  \frac{\sin(\xi/2)}{\xi/2} \right)^{2N}
\prod_{j=1}^{\infty}\ml(2^{-j}\xi).
$$
By (\ref{decayL}), we have
 \bea
 |\hat{\phi}(\xi)|
 &\leq&
 \left( \frac{2}{\xi} \right)^{2N}\exp(C_2/3) q_1q_2^{J-1}
 |\xi|^{\log_2(q_1^{1/(J-1)}q_2)}\\
 &=& C_3 |\xi|^{-2N+\log_2(q_1^{1/(J-1)}q_2)}.
 \eea
 where $C_3 = 4^N\exp(C_2/3) q_1q_2^{J-1}$.
\end{proof}

\begin{cor}
For positive integers $N,l\in\mathbb N$  with $l<N$, let $\phi$ be
the pseudo splines of Type II with order $(N,l)$. Then
 \be\label{philow}
 |\hat{\phi}(\xi)|\leq C \min\{1,|\xi|^{-2N+\kappa}\},
 \ee
where $\kappa = \log(P(3/4))/\log2$.
\end{cor}
\begin{proof}
The proof is a straightforward consequence of Theorem \ref{phiup}.
\end{proof}
\begin{rem}
The above results were first proved by Dong and Shen in \cite{DS1}.
The decay of the Fourier transform is optimal. We give estimation
(\ref{phic}) with explicit constant which is important for the frame
upbound of the shearlet frame. Table 1 gives the decay rate of the
Fourier transform of pseudo splines of Type II with order $(N,l)$
for $2\leq N\leq 9$ and $0\leq l\leq N$. The decay rate  of the
Fourier transform of pseudo splines of Type II with order $(N,l)$ is
$\beta_{N,l}/2$.
\end{rem}
\begin{table}
\caption{Decay rate $\beta_{N,l}=2N-\kappa$ of pseudo splines of
type II for $2\leq N\leq 9$}
\begin{center}
\begin{tabular}{llllllllll} \hline
 $(N,l)$&$N=2$  &  $N=3$   &  $N=4$   &  $N=5$    &  $N=6$     & $N=7$      & $N=8$       & $N=9$     &  \\ \hline
 $l=0$  & $4.00000$ & $6.00000$& $8.00000$& $10.0000$&$12.0000$ & $14.0000$ & $16.0000$  & $ 18.0000$&        \\
 $l=1$  & $2.67807$  & $4.29956$& $6.00000$& $7.75207 $&  $9.54057$ & $11.3561$ & $13.1927$  & $ 15.0458$&        \\ 
 $l=2$  &  & $3.27208$& $4.73321$& $6.27890 $&  $7.88626$ & $9.54057$  & $11.2318$  & $ 12.9530$&        \\ 
 $l=3$  &  &          & $3.82507$& $5.19506 $&  $6.64465$ & $8.15608$  & $9.71691$   & $11.3181$ &        \\ 
 $l=4$  &  &          &          & $4.35316 $&  $5.66363$ & $7.04717$  & $8.48992$   & $9.98156$ &        \\ 
 $l=5$  &  &          &          &           &  $4.86449$ & $6.13261$  & $7.46770$   & $8.85865$ &        \\ 
 $l=6$  &  &          &          &           &            & $5.36349$  & $6.59988$   & $7.89780$  &        \\ 
 $l=7$  &  &          &          &           &            &            & $5.85310$   & $7.06473$ &        \\ 
 $l=8$  &  &          &          &           &            &            &             & $6.33529$ &        \\ \hline
\end{tabular}
\end{center}
\end{table}

\subsection{Shearlet Frames}
In the last subsection, we construct the shearlet frames based on
pseudo splines. We give a weaker condition for constructing shearlet
frames than Theorem 4.9 in \cite{KGL2}.

For function $\phi, \psi, \tp\in L^2(\R)$, we define $\Theta: \R
\times \R \rightarrow \mathbb{R}$ by
$$
\Theta(\xi,\omega) = |\hat{\phi}(\xi)||\hat{\phi}(\xi+\omega)|+
\Theta_1(\xi,\omega) + \Theta_2(\xi,\omega),
$$
where
$$
\Theta_1(\xi,\omega) = \sum_{j\geq 0} \sum_{|k|\leq \lceil 2^{j/2}
\rceil} |\hat{\psi}(S^T_k A_{2^{-j}}\xi)||\hat{\psi}(S^T_k
A_{2^{-j}}\xi+\omega)|
$$
and
$$
\Theta_2(\xi,\omega) = \sum_{j\geq 0} \sum_{|k|\leq \lceil 2^{j/2}
\rceil} |\hat{\tp}(S_k \ta_{2^{-j}}\xi)||\hat{\tp}(S_k
\ta_{2^{-j}}\xi+\omega)|.
$$
The  sufficient conditions for the construction of shearlet frames
were given by \cite{KGL2}.
\begin{lem}\label{thm1}[Theorem 3.4 in \cite{KGL2}]
let $\phi,\psi\in L^2(\R)$ be functions such that
$$
 \hat{\phi}(\xi_1,\xi_2) \leq C_2 \min\{1,|\xi_1|^{-\gamma}\}\min\{1, |\xi_2|^{-\gamma}\}
$$
and
$$
 |\hat{\psi}(\xi_1,\xi_2)| \leq C_1 \min\{1, |\xi_1|^{\alpha}\}
 \min\{1,|\xi_1|^{-\gamma}\}\min\{1, |\xi_2|^{-\gamma}\},
$$
for some positive constants $C_1,C_2<\infty$ and $\alpha>\gamma>3$.
Define $\tp(x_1,x_2)=\psi(x_2,x_1)$ and let $L_{\inf}$ be defined by
$ L_{\inf} = \rm{ess}\inf_{\xi\in\R} \Theta(\xi,0)$. Suppose that
there is a constant $\tilde{L}_{\inf}>0$ such that
$0<\tilde{L}_{\inf}\leq L_{\inf}$. Then there exist a sampling
parameter $c=(c_1,c_2)$ with $c_1=c_2$ such that
$SH(\phi,\psi,\tilde{\psi};c)$ forms a frame for $L^2(\R)$
\end{lem}


Now we state the main results in this section.

\begin{thm}\label{thmframe}
For $0<\alpha<\pi/2$. Let $\widehat{a_{N_1,l_1}}(\xi),$
$\widehat{a_{N_2,l_2}}(\xi)$  be the mask of the pseudo splines of
type II. For $l_2=0$, let $N_1>N_2>2$. For $l_2>0$, let $N_1 \geq
N_2 >2$. Let $\hat{\phi}(\xi)$ be the associated refinable function
with $\widehat{a_{N_2,l_2}}(\xi)$. Define the shearlet $\psi\in
L^2(\R)$ by
  \be\label{psiconstruction}
  \hat{\psi}(2\xi) =
  \hat{b}(\xi_1)\hat{\phi}(\xi_1)\hat{\phi}(\xi_2),\quad \xi =
  (\xi_1,\xi_2)\in \R,
  \ee
where  $\hat{b}(\xi_1) =
e^{-i\xi_1}\overline{\widehat{a_{N_1,l_1}}(\xi_1+\pi)}$. Then for
given $\alpha\in (0, \pi/2)$, there exits a sampling constant
$\hat{c}>0$ such that the shearlet system $\Psi(\psi;c)$ forms a
frame for
 $
 \{ f\in L_2(\R) :
 \mathrm{supp} \hat{f}\in \mathcal{C}_1(\alpha)\cup\mathcal{C}_3(\alpha) \}
 $
with $c=(c_1,c_2)\in (\mathbb R_+)^2$ and $c_2\leq c_1\leq \hat{c}$.
\end{thm}
\begin{proof}
We first give the lower bound of the shearlet $\hat{\psi}(\xi)$. Let
$0<\alpha<\pi/2$. By Theorem \ref{thmphilow}, we obtain
 \be\label{psi1}
 |\hat{\phi}(\xi_1)|\geq C'_4 \chi_{[-2\alpha,2\alpha]}(\xi_1)
  \quad\text{and}\quad
 |\hat{\phi}(\xi_2)|\geq C_4 \chi_{[-\alpha,\alpha]}(\xi_2).
 \ee
From (\ref{blowup}), we have
 \be\label{psi2}
|\hat{b}(\xi_1)|\geq|\hat{b}(\alpha)| \chi_{[-2\alpha,
-\alpha]\cup[\alpha,2\alpha]}(\xi_1).
 \ee
Combine (\ref{psi1}) and (\ref{psi2}), we obtain
 $$
 |\hat{\psi}(\xi)|  =  |\hat{b}(\xi_1) \hat{\phi}(\xi_1)\hat{\phi}(\xi_2)|
 \geq |\hat{b}(\alpha)|C'_4C_4\cdot \chi_{\Omega}(\xi)
 $$
where $\Omega = \{ \xi=(\xi_1,\xi_2), \quad \xi_1 \in
[-2\alpha,-\alpha]\cup[\alpha,2\alpha],\quad \xi_2 \in
[-\alpha,\alpha] \}.$ Since
$$
\bigcup_{j=0}^{\infty}\bigcup_{|k|\leq \lceil 2^{j/2} \rceil }
A_{2^j}S^T_k \Omega = \mathcal{C},
$$
we conclude that
$$
\Phi(\xi,0)=\sum_{j,k}|\hat{\psi}(S^T_kA_{2^{-j}}\xi)|^2 \geq
\tilde{L}_{\inf} \chi_{\Omega}(S_k^TA_{2^{-j}}\xi)\geq
\tilde{L}_{\inf}\quad \text{on}\quad \mathcal{C},
$$
where $\tilde{L}_{\inf} = (|\hat{b}(\alpha)|C'_4C_4)^2$.

By (\ref{blowup}) and  (\ref{philow}), we obtain the upper bound of
the shearlet $\hat{\psi}(\xi)$
 \be\label{psiup}
|\hat{\psi}(\xi)| \leq C \min\{1,|\xi_1|^{2N_1}\}\cdot
 \min\{1,|\xi_1|^{-2N_2+\kappa}\}\cdot
 \min\{1,|\xi_2|^{-2N_2+\kappa}\}.
 \ee
Let $\beta_{N,l}=2N-\kappa$, for fixed $N$, $\beta_{N,l}$ decreases
as $l$ increases. For fixed $l$, $\beta_{N,l}$ increases as $N$
increases. When $l=N-1$, $\beta_{N,l}$ increases as $N$ increase
\cite{DS1}. From the table 1, it easy to see that
$$
2N_1 > 2N_2 >3, \quad \text{for} \quad l_2=0
$$
and
$$
2N_1 > 2N_2-\kappa >3, \quad \text{for} \quad l_2>0.
$$
Hence, the condition of Theorem \ref{thm1} holds. Te theorem is
proved.
\end{proof}
\begin{rem}
(\ref{psiconstruction}) is a standard construction for separable
wavelet basis of $L_2(\mR^2)$ \cite{DTen}. Therefore, Theorem
\ref{thmframe} gives a connection between the shearlet analysis and
wavelet analysis and also, we hope, enriches the theory of shearlet
analysis.
\end{rem}

Since the decay rate of $|_1\hat{\phi}(\xi)|$ is half of that of
$|_2\hat{\phi}(\xi)|$. The results on pseudo splines of type I
follow directly from Theorem \ref{thmframe}.

\begin{thm}
For $0<\alpha<\pi/2$. Let $\widehat{a_{N_1,l_1}}(\xi)$,
$\widehat{a_{N_2,l_2}}(\xi)$ be the mask of the pseudo splines of
type I. For $l_2=0$, let $N_1>N_2>3$. For $l_2>0$, let $N_1 \geq N_2
>8$. Let $\hat{\phi}(\xi)$ be the associated refinable function with
$\widehat{a_{N_2,l_2}}(\xi)$. Define a shearlet $\psi\in L^2(\R)$ by
  $$
  \hat{\psi}(2\xi) =
  \hat{b}(\xi_1)\hat{\phi}(\xi_1)\hat{\phi}(\xi_2),\quad \xi =
  (\xi_1,\xi_2)\in \R,
  $$
where  $\hat{b}(\xi_1) =
e^{-i\xi_1}\overline{\widehat{a_{N_1,l_1}}(\xi_1+\pi)}$. Then there
exits a sampling constant $\hat{c}>0$ such that the shearlet system
$\Psi(\psi;c)$ forms a frame for
 $
 \{ f\in L_2(\R) :
 \mathrm{supp} \hat{f}\in \mathcal{C}(\alpha)_1\cup\mathcal{C}(\alpha)_3  \}
 $
for any sampling matrix $M_c$ with $c=(c_1,c_2)\in (\mathbb R_+)^2$
and $c_2\leq c_1\leq \hat{c}$.
\end{thm}

\begin{rem}
In \cite{KGL2}, P. Kittipoom, G. Kutyniok, and W.-Q Lim gave a
similar construction of the shearlet from pseudo splines of type I
with order $(N,l)$ such that $l>10$ and $\frac{3l}{2}\leq N\leq
3l-2$. Compared with their results, our construction have smaller
support. Moreover, we prove for any given $0<\alpha<\pi/2$, one can
construct shearlet frame. This gives more feasible choice in
applications.
\end{rem}

\section{Examples}
To illustrate our results, we give two examples in this section.
\begin{exa}\label{exa1}
Let $B_3$ be the B-spline function of order $3$.  Define the
shearlet
$$
\hat{\psi}(2\xi)=2^{-4}e^{-i\xi_1} (1-e^{-i\xi_1})^{4}
\widehat{B_{3}}(\xi_1)\widehat{B_{3}}(\xi_2)\quad \xi =
  (\xi_1,\xi_2)\in \R .
$$
Then for any given $0<\alpha<\pi/2$, there exits a sampling constant
$\hat{c}>0$ such that the shearlet system $\Psi(\psi;c)$ forms a
frame for
 $
 \{ f\in L_2(\R) :
 \mathrm{supp} \hat{f}\in \mathcal{C}_1(\alpha)\cup\mathcal{C}_3(\alpha)  \}
 $
for  $c_2\leq c_1\leq \hat{c}$.
\begin{figure}
\centering{}
 \includegraphics[height=8cm]{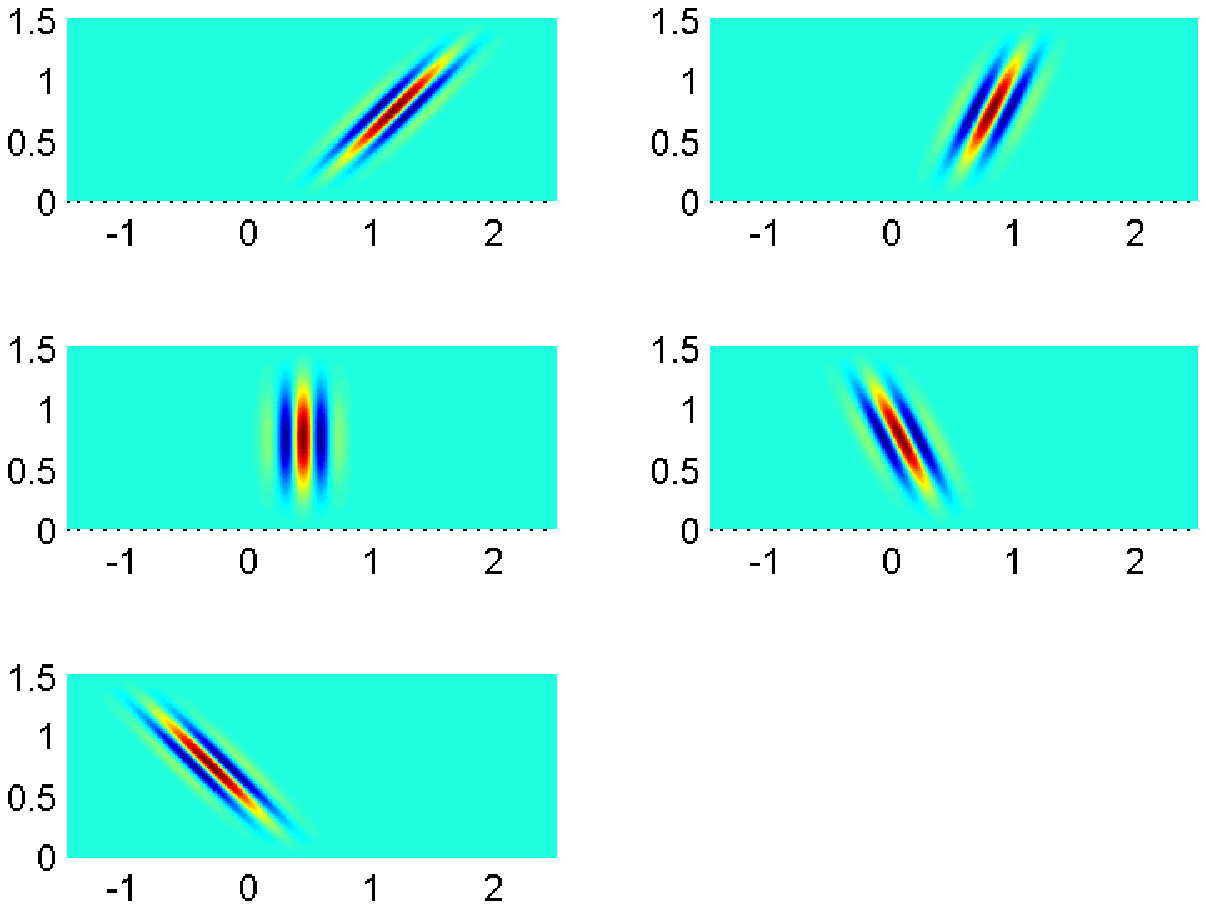}
 \caption{$\psi(S_k A_{2^2}\cdot x)$, $k=-2,-1,0,1,2$, in Example \ref{exa1}}
 \label{fig2}
\end{figure}
\end{exa}


\begin{exa}\label{exa2}
Let $B_4$ be the B-spline function of order $4$.  Define the
shearlet
$$
\widehat{\psi_1}(2\xi)=2^{-6}e^{-i\xi_1} (1-e^{-i\xi_1})^{6}
\widehat{B_{4}}(\xi_1)\widehat{B_{4}}(\xi_2), \quad \xi =
  (\xi_1,\xi_2)\in \R.
$$
Let $\phi(x) = B_{4}(x_1)B_{4}(x_2)$ and
$\psi_2(x_1,x_2)=\psi_1(x_1,x_2)$. Then there exit  sampling
constant $c>0$ such that the shearlet system $\Psi(\phi,
\psi_1,\psi_2;c)$ provides (almost) optimally sparse approximations
of function $f\in\mathscr{E}^2(v)$
\begin{figure}
\centering{}
 \includegraphics[height=8cm]{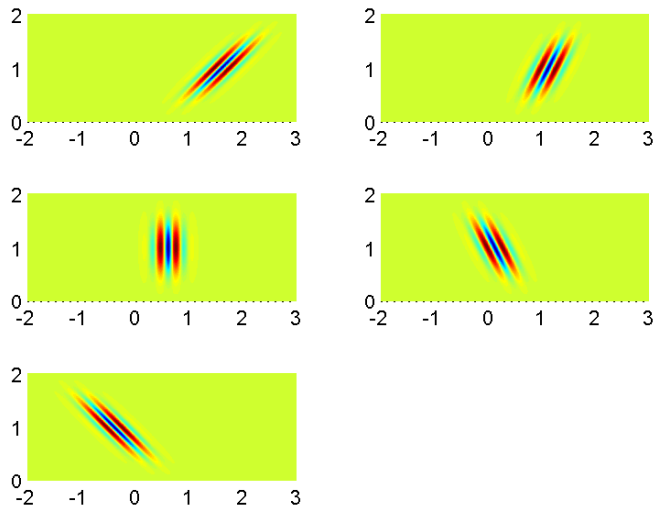}
 \caption{$\psi(S_k A_{2^2}\cdot x)$, $k=-2,-1,0,1,2$, in Example \ref{exa2}}
 \label{fig2}
\end{figure}
\end{exa}

\begin{ac}
This work is supported  by NSF of China under grant numbers
10771190, 10971189 and Zhejiang Provincial NSF of China under grants
number Y6090091.
\end{ac}


\end{document}